	\documentclass[12pt, a4paper]{article}
	
 	\usepackage{latexsym}
	\usepackage{amsmath}	
	\usepackage{epsfig}	
	\usepackage{times}	
	\usepackage{amssymb}
	
 	\newtheorem{thr}{Theorem}[section]
 	\newtheorem{prop}{Proposition}[section]

\begin{document}
 	
 	\centerline{\Large{\bf Dimension theoretical study on skew product maps with}}
 	\centerline{}
 	\centerline{\Large{\bf coupled-expanding property}}
 	\centerline{}
 	\centerline{}
 	\centerline{\textsuperscript{1}Jinhyon Kim, \textsuperscript{2}Hyonhui Ju}
 	\centerline{}
 	{\small \centerline{{ Department of Mathematics,  \textbf{Kim Il Sung} University, Pyongyang, D.P.R.Korea}}}
	{\small \centerline{{E-mail address:  jhyonhui@yahoo.com}}}
	
	\centerline{}
	\centerline{}
	\begin{abstract}
	We discuss on some families of skew product maps on a square. 
	For a kind of skew product maps with coupled-expanding property, we estimate Hausdorff dimension of its attractor. 
	And we prove that there exists an ergodic measure with full Hausdorff dimension for this kind of skew product maps.
	\end{abstract}

	\section{Introduction}
	In this paper we study dimensional theoretical properties of a kind of skew product maps with coupled-expanding property.
	By dimensional theoretical properties we mean Hausdorff dimension of invariant sets of dynamical systems and 
	ergodic measures on these sets. 

	In the modern theory of dynamical systems geometrical invariants like 
	Hausdorff dimension of invariant sets and measures seem to have their place beside classical invariants like 
	entropy and Lyapnov exponents \cite{ne2}. In the dimension theoretical study on dynamical systems, following two 
	problems are regarded as important. The first one is to estimate Hausdorff dimension of invariant sets. 
	The second one is the problem about existence of ergodic measure with full Hausdorff dimension in these invariant sets. 
	If such an ergodic measure does not exist, it is problem whether at least the variational principle for Hausdorff dimension holds or not. 
	The variational principle means that there exists a sequence of ergodic measures such that the Hausdorff dimension 
	of the invariant set can be approximated by their Hausdorff dimensions. 
	It seems to be well accepted by experts that these questions are of great important in developing the dimension theory of 
	dynamical systems \cite{pes}.
	
	The problem of estimating the Hausdorff dimension of invariant subsets has been studied in several papers for different 
	classes of dynamical system (a good overview can be found in \cite{hof, pes}). 
	In \cite{hof}, for some noninjective skew product maps defined on a square without finite Markov partition, 
	the author estimated Hausdorff dimension of attractor under some conditions. 
	In generally the skew product map $F: Q \to Q$ with $Q=[a,b] \times [0,1]$ is defined by $F(x,y)=(T(x),g(x,y))$, 
	where $F: [a,b] \to [a,b]$ is piecewise monotonic and the map $y \mapsto g(x,y)$ is a contraction on $[0,1]$ for all $x \in [a,b]$. 
	In \cite{hof}, it was proved that for the skew product map $F(x,y)=\Big(T(x), \varphi (x)+\lambda \big( y-\frac{1}{2} \big) \Big)$, 
	Hausdorff dimension of the attractor $\Delta:= \bigcap_{n=0}^{\infty}F^n(Q)$ of $F$ is 
	$\textnormal{dim}_H(\Lambda)=1+\frac{h_{top}(T)}{-\textnormal{log} \lambda}$ under some conditions, 
	where $0<\lambda <\textnormal{min} \big\{ \frac{1}{2}, \textnormal{inf} \frac{1}{|T|} \big\}$ and $h_{top}(T)$ 
	is the topological entropy of $T$. 
	
	On the other hand Baker's map which can be considered as a special form of the skew product map is well known as a simple 
	example of a "chaotic" dynamical system. When it comes to dimensional theoretical study for this map, in \cite{ne2} 
	it was considered that, for parameters in some domain, Baker's transformation has ergodic measures with full 
	Hausdorff dimension in its invariant sets but, for parameters in some other domain, variational principle for 
	Hausdorff dimension doesn't hold, that is, there is not a sequence of ergodic measures of which Hausdorff dimension 
	approximates the Hausdorff dimension of the invariant set.
	In \cite{ne1} more improved results than one in \cite{ne2} on the dimensional theoretical 
	properties of the generalized Baker's map were obtained.
	
	According to our search on prior works, it seems that there is no result considered for the existence problem of ergodic 
	measure with full Hausdorff dimension or for the problem on variational principle for Hausdorff dimension for the 
	skew product map yet, whereas such work has done for the Baker's map. And in \cite{hof} rather 
	strong condition including continuously differentiability for the map $T$ and $\varphi$  was supposed. 
	When it comes to coupled-expanding property, strictly coupled-expanding map implies "chaos" in several senses (see \cite{ric, shi}). 

	From the above mentioned results, in this paper we study on dimensional theoretical properties for a kind of skew product maps 
	$F(x,y)=\Big(T(x), \varphi (x)+\lambda \big( y-\frac{1}{2} \big) \Big)$ where $T$ is a coupled-expanding map. 
	We estimate Hausdorff dimension of the attractor of the skew product map $F$ (in section 3) and prove that there exists 
	ergodic measure with full Hausdorff dimension in the attractor for the map $F$ (in section 4).
	\section{Basic notation}
	Let $X$ be a metric space. Recall that Hausdorff dimension of $A \subset X$ is given by 
	\begin{equation*}
	\textnormal{dim}_H(A) = \textnormal{inf} \{s: H^s(A)=0 \} = \textnormal{sup} \{s: H^s(A)=\infty \},
	\end{equation*}
	where $H^s(A)$ denotes the $s$-dimensional Hausdorff measure, i.e.
	\begin{equation*}
	H^s(A)=\lim_{\varepsilon \to 0}\textnormal{inf} \Big\{\sum \vert U_i \vert ^s: A \subset \bigcup U_i,  \vert U_i \vert < \varepsilon \Big\},
	\end{equation*}
	and $\vert U_i \vert$ denotes the diameter of a covering element $U_i$ (see \cite{ne3}).
	
	Let $N$ be a positive integer and $S_i: X \to X (1 \leq i \leq N)$ be a contraction similitude with contraction ratio
	$c_i (0<c_i<1)$, i.e. $d(S_i(x), S_i(y))=c_id(x,y)$ for any $x,y \in X$. A set $K \subset X$ is said to be invariant set 
	with respect to $S:=\{S_1, \cdots, S_N \}$ if $K=\sum_{i=1}^N S_i(K)$. If there exists a finite set $S=\{S_1, \cdots, S_N \}$ 
	of contraction similitudes with contraction ratio $c_i (0<c_i<1)$ such that $K$ is a invariant set with respect to 
	$S=\{S_1, \cdots, S_N \}$, and if $H^{\alpha}(K)>0, H^{\alpha}(K_i \cap K_j)=0 (i \neq j)$ for $\alpha >0$ 
	such that $\sum_{i=1}^N c_i^{\alpha}=1$, then $K$ is said to be self-similar set, where $H^{\alpha}$ is 
	$\alpha$-dimensional Hausdorff measure and $K_i=S_i(K)$ (see \cite{hut}).
	Let  $D \subset X$ and $f: D \subset X \to X$. If there exists $m \geq 2$ subsets $V_i (1 \leq i \leq m )$  of $D$ 
	with $V_i \cap V_j = \emptyset$ and $d(V_i, V_j) >0$  for all $1 \leq i \neq j \leq m$, such that 
	\begin{equation*}
	f(V_i) \supset \bigcup_{j=1}^m V_j, \quad 1 \leq i \leq m,
	\end{equation*}
	then $f$ is said to be coupled-expanding map in $V_i, 1 \leq i \leq m$ (see \cite{shi}).
	\section{Hausdorff dimension of attractor of a skew product map with coupled expanding property}

	\begin{thr}
	Let $[a,b] \subset \mathbb{R}$ be an interval. Suppose that  $T: [a,b] \to [a,b]$ is monotonic on 
	$[a,c]\cup [d,b] (a<c<d<b)$ and satisfies $T([a,c])=T([d,b])=[a,b]$ and $T((a,c)) \subset \{a,b \}$. 
	Assume that $0< \lambda <1/2$ and that $\varphi : [a,b] \to [1/2, 1- \lambda /2]$ is strictly monotonic 
	and satisfies $\vert \varphi (c) - \varphi (d) \vert > \lambda$. Then the attractor 
	$\Lambda :=\bigcap_{n=0}^{\infty}F^n(Q)$ of the map $F$ defined by $F: Q \to Q, 
	F(x,y)=\Big(T(x), \varphi (x)+\lambda \big( y-\frac{1}{2} \big) \Big), Q=[a,b] \times [0,1]$, has Hausdorff dimension
	\begin{equation}
	\textnormal{dim}_H (\Lambda) = 1+\frac{\textnormal{log} 2}{-\textnormal{log} \lambda}.
	\end{equation}
 	\end{thr}
	
	\noindent \textbf{(Remark 1)} From the conditions of Theorem 3.1, the map $T: [a,b] \to [a,b]$ is strictly expanding map 
	(see \cite{shi}) and therefore the entropy of $T$ is  $\textnormal{log}2$ (\cite{ric}). So we can see that above expression (1) is 
	analogue to the one of \cite{hof}.\\

%%%%%%%%%%%%%    proof of Theorem 3.1  %%%%%%%%%%

	\noindent \textit{Proof of Theorem 3.1.}  For convenience of expressing, define follows;
	\begin{equation*}
	V_1:=[a,c], ~~ V_2:=[d,b], ~~ \Sigma_2^+:=\{1,2\}^{\mathbb{N} \cup \{ 0 \}} ~~ \Delta_i:=V_i \times [0,1] (i=1, 2).
	\end{equation*}
 	And set 
	\begin{equation*}
	S_{i_0 \cdots i_{p-1}}:=F^p (\Delta_{i_{p-1}} \cap \cdots \cap F^{-(p-1)} \Delta_{i_0}) ~~ (i_j \in \{ 1,2 \}, ~ j \in \{ 0, \cdots , p-1 \}).
	\end{equation*}
	Then, for any $p \in \mathbb{N}$, it follows that $S_{i_0 \cdots i_{p}} \subset S_{i_0 \cdots i_{p-1}}$ since
	\begin{eqnarray*}
	S_{i_0 \cdots i_{p}} & =  &F^{p+1} (\Delta_{i_{p}} \cap \cdots \cap F^{-p} \Delta_{i_0}) \subset \\
					& \subset & F^p (F\Delta_{i_{p}} \cap \cdots \cap F^{-(p-1)} \Delta_{i_0}) \subset S_{i_0 \cdots i_{p-1}}.
	\end{eqnarray*}  
	On the other hand, we have 
	\begin{equation*}
	S_{i_0 \cdots i_{p-1}}=F^p \big( (V_{i_{p-1}} \cap \cdots \cap T^{-(p-1)} V_{i_0}) \times [0,1] \big)
	\end{equation*}
 	from
	\begin{equation*}
	\Delta_{i_{p-1}} \cap \cdots \cap F^{-(p-1)} \Delta_{i_0} = \big( V_{i_{p-1}} \cap \cdots \cap T^{-(p-1)} V_{i_0} \big) \times [0,1].
	\end{equation*}
	From the definition of the map $F$, it follows that
 	\begin{equation*}
	F^p (x,y)=\big( T^p x, \varphi(T^{p-1}x)+ \lambda \varphi (T^{p-2}x)+ \cdots +\lambda^{p-1} \varphi (x) + \lambda ^p y- \frac{1}{2}(\lambda^p+ \cdots +\lambda) \big).
	\end{equation*}
	Since $T^p \big( V_{i_{p-1}} \cap \cdots \cap T^{-(p-1)} V_{i_0} \big)=[a,b]$, the set $S_{i_0 \cdots i_{p-1}}$ forms a band 
	with bandwidth $\lambda^p$ across the lines $x=a$ and  $x=b$. Therefore $S_{i_0 \cdots i_{p} \cdots}$  is a piece of the curve 
	across the lines $x=a$ and  $x=b$. Then the set $\Lambda$ consists of uncountable pieces of curve, that is, 
	\begin{equation*}
	\Lambda = \bigcup \big\{ S_{i_0 \cdots i_{p} \cdots} : I=i_0 \cdots i_{p} \cdots \in \Sigma_2^+ \big\}
	\end{equation*}
 	Next, let’s consider the skew product map $F_u$ defined by 
	\begin{equation*}
	F_u (x,y)=\Big( T(x), q_i+ \lambda \big( y-\frac{1}{2} \big) \Big), ~~ x \in V_i,
	\end{equation*}
	where $q_i$ is defined by $q_1=\varphi (a)$ and $q_2=\varphi (b)$.

	With the same way as we have considered for $F$ above, we can consider the construction of the attractor of 
	$F_u: Q \to Q$ by $\Lambda_u = \bigcap_{n=0}^{\infty} F_u^n Q$.
	For $I \in \Sigma_2^+$ let $S_I^u$ be the set for $F_u$ corresponding to $S_I$. Then we can see that the set $S_I^u$ 
	is a segment with length $b-a$.
	
	Now let's consider the map $\Phi_u: \Lambda_u \to \Lambda$ defined by $(t,y) \in S_I^u \mapsto (t, y') \in S_I$. 
	The map $\Phi_u$ maps the point at which $S_I^u$  intersects the line $x=t$, to the point at which $S_I$  intersects the 
	line $x=t$, where $t \in [a,b]$. Then $\Phi$  is Lipschitz towards $y$-ax direction, that is, for any fixing $x \in (a,b)$  we have 
	\begin{equation*}
	\vert \Phi_u (x,y)-\Phi_u (x,y') \vert \leq \vert y-y' \vert.
	\end{equation*}
	Therefore 
	\begin{equation*}
	\textnormal{dim}_H \Lambda \leq \textnormal{dim}_H \Lambda_u
	\end{equation*}
 	by using  the result of \cite{ste} (see Fig.1).

	\begin{figure}[tbhp]
	\centering
 	\includegraphics[width=85mm]{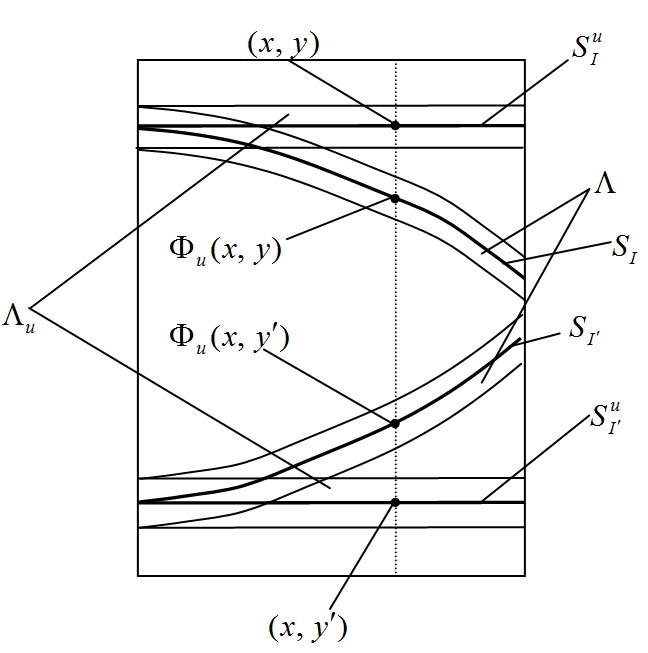}\\ 
  	\end{figure}

	\centerline{{\bf Fig. 1}. Imagine of the point $(x,y) \in \Lambda_u$ by $\Phi_u$}

	\vspace{0.5cm}

	Let $C_u \subset [0,1]$ be the intersection of $\Lambda_u$ and line $x=t$. Then $C_u$  is self-similar set with 
	respect to the contraction similitudes with contraction ratio $\lambda$  and we have $\Lambda_u=[a,b] \times C_u$ (see \cite{hut}).
	Since $\textnormal{dim}_H C_u = \frac{\textnormal{log} 2}{-\textnormal{log} \lambda}$, we have 
	\begin{equation*}
	\textnormal{dim}_H \Lambda_u = 1+ \frac{\textnormal{log} 2}{-\textnormal{log} \lambda}.
	\end{equation*}
	This implies that 
	\begin{equation*}
	\textnormal{dim}_H \Lambda \leq 1+ \frac{\textnormal{log} 2}{-\textnormal{log} \lambda}.
	\end{equation*}
	To prove the opposite inequality, let's consider the map defined by
	\begin{equation*}
	F_l: Q \to Q, ~~ F_l (x,y)=\Big( T(x), \gamma_i + \lambda \big( y- \frac{1}{2} \big) \Big), ~~ x \in V_i,
	\end{equation*}
	where $\gamma_1 = \varphi (c)$ and $\gamma_2 = \varphi (d)$. 

	By the same procedure as we just have done for $F$ and $F_u$, we can consider the construction of the attractor 
	$\Lambda_l$ for $F_l$.  Let $S_I^l$ be the set corresponding to $S_I$, then $S_I^l$  is also a segment with length $b-a$.
	Let $\Phi_l: \Lambda \to \Lambda_l$ be a map which maps intersecting point of the line $x=t$ and $S_I$  to 
	intersecting point of the line $x=t$ and $S_I^l$. Then $\Phi_l$ is also Lipshitz towards $y$-ax direction 
	like as $\Phi_u$ (see Fig.2). 
	
	\begin{figure}[tbhp]
	\centering
 	\includegraphics[width=85mm]{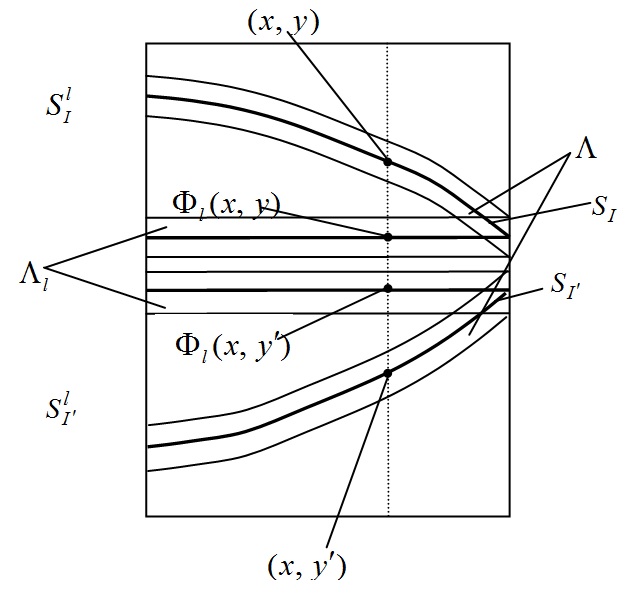}\\ 
  	\end{figure}

	\centerline{{\bf Fig. 2}. Imagine of the point $(x,y) \in \Lambda$ by $\Phi_l$}

	\vspace{0.5cm}

	\noindent Therefore 
	\begin{equation*}
	\textnormal{dim}_H \Lambda_l \leq \textnormal{dim}_H \Lambda.
	\end{equation*}
	Let  $C_l$ be the intersection of the line $x=t$ and $\Lambda_l$. Then $C_l$ is also self-similar set with respect to 
	the contraction similitudes with contraction ratio $\lambda$. This means that 
	$\textnormal{dim}_H C_l = \frac{\textnormal{log} 2}{-\textnormal{log} \lambda}$. Since $\Lambda_l = [a,b] \times C_l$, 
	we have 
	\begin{equation}
	\textnormal{dim}_H \Lambda \geq 1+ \frac{\textnormal{log} 2}{-\textnormal{log} \lambda}.
	\end{equation}
	From (1) and (2) it follows that 
	\begin{equation*}
	\textnormal{dim}_H \Lambda = 1+ \frac{\textnormal{log} 2}{-\textnormal{log} \lambda}.
	\end{equation*}
	The proof is thus complete.
	\section{The existence of ergodic measure with full Hausdorff dimension on attractor}
	Let $M(X,f)$ be the space of all $f$-invariant and ergodic measures on $X$.
	\begin{prop}
	Suppose that the dynamical system $(\Lambda_2, g)$ is topologically semi-conjugate to $(\Lambda_1, f)$, that is, 
	there is a surjection $\varphi: \Lambda_2 \to \Lambda_1$ continuous such that $f \circ \varphi = \varphi \circ g$. 
	Then, If $\mu$ is $g$-invariant and ergodic measure on $\Lambda_2$, then $\mu \circ \varphi^{-1}$  is $f$-invariant 
	and ergodic measure on $\Lambda_1$.
	\end{prop}
	\noindent \textit{Proof}. Assume that $f^{-1}(E)=E, E \subset \Lambda_1$. 
	Then $\varphi^{-1} \big( f^{-1}(E) \big) = \varphi^{-1} (E)$  and $g^{-1} \big(\varphi^{-1}(E) \big) = \varphi^{-1} (E)$  
	because $f \circ \varphi = \varphi \circ g$. Since $\mu$ is $g$-ergodic, $\mu \big( \varphi^{-1}(E) \big) = 0$ or $1$. 
	Therefore $\mu \circ \varphi^{-1}$ is $f$-ergodic measure on $\Lambda_1$. And let $A \subset \Lambda_1$ be a 
	measurable set. Since $\mu$ is $g$-invariant, we have  
	\begin{eqnarray*}
	\mu \circ \varphi^{-1} \big( f^{-1}(A) \big) & = & \mu \big( \varphi^{-1} \big( f^{-1}(A) \big) \big) = \mu \big(g^{-1} \big( \varphi^{-1}(A) \big) \big) \\
			& = & \mu \big( \varphi^{-1}(A) \big)  = \mu \circ \varphi^{-1}(A).
	\end{eqnarray*}
	Therefore the proof is complete.\\

	\begin{thr}
	Suppose that the map  
	\begin{equation*}
	F: Q \to Q, ~~ F(x,y)=\Big( T(x), \varphi (x)+\lambda \big( y- \frac{1}{2} \big) \Big)
	\end{equation*}
	satisfies the assumptions of Theorem 3.1 and that $T$ is linear on intervals $[a, c]$ and $[d,b]$. 
	Then there exists a $F$-invariant ergodic measure with full dimension on the attractor $\Lambda = \bigcap_{n=0}^{\infty} F^nQ$.
 	\end{thr}
	\noindent \textit{Proof}. Let  $L$  be a normalized  Lebesgue measure on interval $[a,b]$. Since $C_l$ is self similar set 
	with respect to two self similar contractions with contraction ratio $\lambda$, it is homeomorphism to $\Sigma_2^+$. 

	Let $\beta$ be a Bernoulli measure on $\Sigma_2^+$  and $\mu_{\beta}$ be a measure corresponding to $\beta$  on $C_l$. 
	Define $\mu_l$ as follows;  
	\begin{equation*}
	\mu_l = L \times \mu_{\beta}.
	\end{equation*}
	Then $\mu_l$  is  $F_l$ ergodic measure defined on $\Lambda_l$  and it is obvious that 
 	\begin{equation*}
	\textnormal{dim}_H \mu_l = 1+ \frac{\textnormal{log} 2}{-\textnormal{log} \lambda}
	\end{equation*}
 	by composing of  $\mu_l$. By using proposition 4.1, it follows that  $\mu:=\mu_l \circ \Phi_l$ is $F$-ergodic measure.
	Let $\mu (E)=1, E \subset \Lambda$. From Theorem 3.1 it follows that
	\begin{equation*}
	\textnormal{dim}_H E \leq \textnormal{dim}_H \Lambda = 1+ \frac{\textnormal{log} 2}{-\textnormal{log} \lambda}.
	\end{equation*}
	Since $\mu_l (\Phi_l E)=1$, from the definition of Hausdorff dimension we have 
	\begin{equation*}
	\textnormal{dim}_H \Phi_l E \geq \textnormal{dim}_H \mu_l = 1+ \frac{\textnormal{log} 2}{-\textnormal{log} \lambda}.
	\end{equation*}
	Remember $\Phi_l$ is Lipshitz towards $y$-ax direction, then we have 
	\begin{equation*}
	\textnormal{dim}_H E \geq \textnormal{dim}_H \Phi_l E \geq 1+ \frac{\textnormal{log} 2}{-\textnormal{log} \lambda}.
	\end{equation*}
	Thus it follows that
	\begin{equation*}
	\textnormal{dim}_H E= 1+ \frac{\textnormal{log} 2}{-\textnormal{log} \lambda}.
	\end{equation*}
	Considering
	\begin{equation*}
	\textnormal{dim}_H \mu = \textnormal{inf} \{ \textnormal{dim}_H E : \mu (E)=1, E \subset \Lambda \},
	\end{equation*}
	we have 
	\begin{equation*}
	\textnormal{dim}_H \mu = 1+ \frac{\textnormal{log} 2}{-\textnormal{log} \lambda} = \textnormal{dim}_H \Lambda.
	\end{equation*}
	The proof is thus complete.
	\section{Conclusion and further study}
	In this paper, for the skew product maps with coupled-expanding property, we estimated Hausdorff dimension 
	of its attractor and proved that there exists an ergodic measure with full Hausdorff dimension for the attractor. 
	In further study we are going to weaken the limit for the skew product map $F$ which is actually assumed 
	to avoid overlapping of its image here. \\

	{\bf Acknowledgement}. We would like to thank anonymous reviewers' help and advice.

	\end{document}